\documentclass[11pt]{article}
\usepackage{graphicx}
\usepackage{amssymb}
\usepackage{amstext}
\usepackage{amsthm}
\usepackage{epstopdf}
\newcommand{\R}{{\mathbb R}}

\newtheorem{theorem}{Theorem}
\newtheorem{corollary}[theorem]{Corollary}
\newtheorem{lemma}[theorem]{Lemma}
\newtheorem{proposition}[theorem]{Proposition}
\newtheorem{remark}{Remark}

\title{On the global attractor of delay differential equations with unimodal feedback}
\author{
 Eduardo Liz
\\
 {\small  Departamento de Matem\'atica Aplicada II,  E.T.S.E. Telecomunicaci\'on }\\
 {\small  Campus Marcosende, Universidade de Vigo, 36310 Vigo,  Spain}\\  ({\small\tt eliz@dma.uvigo.es})
 \medskip \\  {\small and }\medskip\\
 Gergely R\"{o}st \\ {\small Analysis and Stochastics Research Group, Hungarian Academy of Sciences}\\
 {\small Bolyai Institute, University of Szeged, H-6720 Szeged, Aradi v\'ertan\'uk tere 1., Hungary}\\
 ({\small\tt rost@math.u-szeged.hu})
 }
\date{\today}

\begin{document}
\maketitle

\begin{abstract} We give bounds for the global attractor of the delay differential equation $ \dot
x(t)=-\mu x(t)+f(x(t-\tau))$, where $f$ is unimodal and has negative
Schwarzian derivative. If $f$ and $\mu$ satisfy certain condition, then,
regardless of the delay, all solutions enter the domain where $f$
is monotone decreasing and the powerful results for delayed
monotone feedback can be applied to describe the asymptotic
behaviour of solutions. In this situation we determine the sharpest interval that contains
the global attractor for any delay. In the absence of that condition, improving
earlier results, we show that if the delay is sufficiently small,
then all solution enter the domain where $f'$ is negative. Our
theorems then are illustrated by numerical examples using
Nicholson's blowflies equation and the Mackey-Glass equation.

\end{abstract}


\section{Introduction}

This note is motivated by a recent paper by G. R\"ost and J. Wu \cite{rw}  about the so-called  delayed recruitment
model defined by the delay differential equation
\begin{equation}
\label{1}
x'(t)=-\mu x(t)+f(x(t-\tau)),
\end{equation}
where $\mu>0$, $\tau>0$, and $f:[0,\infty)\to [0,\infty)$ is a continuous function.

In particular, they consider the case when $f$ is a unimodal
function, which is the situation for the famous Nicholson's
blowflies equation and the Mackey-Glass model.

In that reference, the authors have proved several results on the
global dynamics of Eq. (\ref{1}), and they also formulated some open
problems. It is our purpose to prove new results in the direction
initiated in \cite{rw}, and also to answer some of the open
questions. We show the applicability of our results for different
cases of the Nicholson's blowflies equation
\begin{equation}
\label{nic}
x'(t)=-\mu x(t)+px(t-\tau)e^{-\gamma x(t-\tau)},
\end{equation}
where $\mu$, $p$, $\gamma$ are positive parameters (see, e.g.,
\cite{ltt} for a biological interpretation); and the Mackey-Glass
equation \cite{mg}
\begin{equation}
\label{mg0}x'(t)=-\mu x(t)+\frac{px(t-\tau)}{1+x(t-\tau)^{n}}.
\end{equation}

Following \cite{rw}, we assume that $f$ is unimodal. More precisely, the following hypothesis will be required:
\begin{enumerate}
\item[(U)] $f(x)\geq 0$ for all $x\geq 0$, $f(0)=0$, and there is a unique $x_0>0$ such that $f'(x)>0$ if
$0\leq x<x_0$, $f'(x_0)=0$, and $f'(x)<0$ if $x>x_0$. Moreover, $f''(x)<0$ if $0\leq x\leq x_0$.
\end{enumerate}

Much is known about the global picture of the dynamics of Eq.
(\ref{1}), when $f$ is a monotone function. However, unimodal
feedback may lead to very complicated and still not completely
understood dynamics. See \cite{kw, kww, lw} and
references thereof.

We shall use the function $g(x)=\mu^{-1}f(x)$. Notice that the
equilibria of (\ref{1}) are the fixed points of $g$, and that, under
condition (U), function $g$ has at most two fixed points $x=0$ and
$x=K>0$.

As in \cite{rw}, we consider nonnegative solutions of (\ref{1}). We recall that for each nonnegative and nonzero
function $\phi\in C=C([-\tau,0],\R)$, there exists a unique solution $x^{\phi}(t)$ of (\ref{1}) such that
$x^{\phi}=\phi$ on $[-\tau,0]$. Moreover,   $x^{\phi}(t)>0\, ,\,\forall\, t>0.$ (See, e.g., \cite[Corollary 12]{gt}.)

It is well known (see, e. g., \cite{gt, rw}) that  all solutions of (\ref{1}) converge to $0$ if $g'(0)\leq 1$,
whereas  the positive equilibrium $K$ is globally attracting for Eq. (\ref{1}) if $g'(K)\geq 0$ (equivalently,
if $K\leq x_0$).

Thus, we shall assume that $g'(0)>1$ and $K>x_0$. Ivanov and
Sharkovsky \cite[Theorem 2.3]{is} proved that an invariant and
attracting interval $[\alpha,\beta]$ for $g$ is also invariant and
attracting for (\ref{1}) for all values of the delay $\tau$, that
is,
$$
\alpha\leq \liminf_{t\to\infty}x(t)\leq \limsup_{t\to\infty}x(t)\leq\beta,
$$
for any nonzero solution $x$ of (\ref{1}). Similar results were proven using a slightly different approach in
\cite{gt, rw}.

It is clear that we can choose $\beta=g(x_0)$,
$\alpha=g(\beta)=g^2(x_0)$ to get an attracting invariant interval
$[\alpha,\beta]$ for the map $g$ (here, and in the following, $g^2$ denotes the composition $g\circ g$). Thus, this interval contains the
global attractor associated to Eq (\ref{1}) for all values of
$\tau$. Using this fact, R\"ost and Wu obtain sufficient conditions
to ensure that every solution of (\ref{1}) enters the domain where
$f'$ is negative. In this case the asymptotic behaviour of the
solutions is governed by monotone delayed feedback and the
comprehensive theory of monotone dynamics is applicable, as it was
demonstrated in \cite{rw}. In particular, since a
Poincar\'e-Bendixson type theorem is available for (\ref{1}) when
$f'$ is negative, this kind of conditions exclude the possibility of
solutions with complicated asymptotic behaviour (the $\omega$-limit
set can only be the positive equilibrium $K$ or a periodic orbit).
We include here the main results of \cite{rw} in this direction
\begin{theorem}
\label{trw}
Every solution of (\ref{1}) enters the domain where $f'$ is negative if any of the following conditions holds:
\begin{enumerate}
\item[{\rm (L)}] $\alpha=g(\beta)=g^2(x_0)> x_0$.
\item[{\rm (L$_{\tau}$)}] $\tau<\tau^{*}:=\displaystyle\frac{\Pi(x_0)-x_0}{\mu(g(x_0)-g^2(x_0))}\; ,$

where $\Pi:(0,K)\to [K,\infty)$ is the inverse of the restriction of $g$ to the interval $[K,\infty)$.
\end{enumerate}
\end{theorem}

Notice that the first condition in Theorem \ref{trw} is independent of the delay, while condition {\rm (L$_{\tau}$)}
shows that even if $f$ is unimodal,  the solutions of (\ref{1}) have the same asymptotic behaviour  as  in the case
of monotone decreasing  feedback for all sufficiently small delay $\tau$.

An open problem suggested in \cite{rw} is the following: under condition (L), find the sharpest invariant and attracting
 interval containing the global attractor of (\ref{1}) for all $\tau$. (Numerical experiments performed in \cite{rw}
 show that $J=[\alpha,\beta]$ seems to be a very sharp bound).
To avoid confusion we remark that the global attractor $\mathcal{A}$ is a subset of the function space
$C([-\tau,0],\R)$, and saying that an interval $[a,b]$ contains the global attractor
we mean that for each $\phi \in \mathcal{A}$, we have
$a \leq \phi(s) \leq b$ for any $s \in [-\tau,0].$

Our main results in this note are the following:
\begin{enumerate}
\item We completely solve this problem in the case of Nicholson's blowflies and Mackey-Glass equations, obtaining an interesting dichotomy result for (\ref{nic}) and (\ref{mg0}) when condition (L) holds (Theorem \ref{tdic}).
\item We give a weaker delay-dependent condition different from  {\rm (L$_{\tau}$)}  under which the statement of
Theorem \ref{trw} remains valid. In other words, we can determine a
$\tau_*$ that is larger than $\tau^*$ in Theorem 1, such that all
solutions enter the domain where $f'$ is negative, if $\tau <
\tau_*$. Moreover, we provide some examples showing that the new
condition significantly improves {\rm (L$_{\tau}$)} in certain
situations.
\end{enumerate}

\section{Main Results}
\setcounter{equation}{0}

In this section, we assume that condition (U) holds, $g'(0)>1$ and $K>x_0$, where $g=\mu^{-1}f$.
Denote
\begin{equation}
\label{sharp}
\bar\alpha=\inf\{A>0\, :\, g^2(A)=A\}\quad ;\quad \bar\beta=\sup\{B>0\, :\, g^2(B)=B\}.
\end{equation}

Since $g'(0)>1$, $\bar\alpha$ and $\bar\beta$ are well defined real numbers.
Assuming that (L) is satisfied, it is clear that $\bar\alpha>x_0$ and hence $g$ is decreasing on
$\bar J:=[\bar\alpha,\bar\beta]$. Moreover, $g(\bar\alpha)=\bar\beta$, and $g(\bar\beta)=\bar\alpha$. As a consequence,
 $\bar J$ is invariant for the map $g$.

\begin{lemma}
\label{l1}
Assume that (L) holds. Then, $\bar J=[\bar\alpha,\bar\beta]$ is an attracting invariant interval for the map $g$.
\end{lemma}
\proof We have already proved that $\bar J$ is invariant. Next, we prove that it is attracting.

Since $h=g^2$ is monotone increasing in $\bar J$, and $\bar \alpha$ and $\bar\beta$ are respectively the minimal and the
maximal fixed points of $h$ in $[\alpha,\beta]$, it follows that
$$
\lim_{n\to\infty}h^n(x)=\bar\alpha,\,\forall\, x\in [\alpha, \bar\alpha]\quad ;\quad
\lim_{n\to\infty}h^n(x)=\bar\beta,\,\forall\, x\in [\bar\beta, \beta].
$$
Since $g([\alpha, \bar\alpha])=[\bar\beta, \beta]$ and $g([\bar\beta, \beta])=[\alpha, \bar\alpha]$, the result follows
from the fact that $[\alpha,\beta]$ is attracting for $g$. \qed

An application of the above mentioned Theorem 2.3 in \cite{is} gives the following result:
\begin{corollary}
\label{c3}
Under condition (L), the interval $\bar J=[\bar\alpha,\bar\beta]$ is an attracting and invariant interval for (\ref{1})
for all values of the delay $\tau$.
\end{corollary}

\begin{remark}
\label{r1}
 If the equilibrium $K$ is globally attracting for $g$, then $\bar\alpha=\bar\beta=K$, and hence $K$ is also a global
 attractor for (\ref{1}) for all values of the delay. This result is \cite[Theorem 2.2]{is}.
\end{remark}

\begin{remark}
\label{r2}
If (L) does not hold, the interval $\bar J$ does not need to be globally attracting for  (\ref{1}).  For example, for
the Nicholson's blowflies equation considered in \cite{rw}
\begin{equation}
\label{nic0}
x'(t)=- 0.05 x(t)+ x(t-\tau)e^{-x(t-\tau)},
\end{equation}
the interval $\bar J$ is given by
$$
\bar J=[\bar\alpha, \bar\beta] \approx [0.4261, 5.5653].
$$
It is easy to check that this interval does not attract every orbit associated to the map $g$, since there is a period
four orbit given by
$$
0.24286\rightarrow 3.80991 \rightarrow 1.6878 \rightarrow 6.24235 \rightarrow 0.24286.
$$
Numerical experiments from \cite{rw} suggest that $\bar J$ also does not attract every orbit of (\ref{nic0}) for large
values of $\tau$.
\end{remark}

We recall here that the nonlinearity $f$ in some important examples of Eq. (\ref{1}) (including the Mackey-Glass and
Nicholson's blowflies models) fulfills the following additional assumption:

\begin{enumerate}
\item[(S)] $f$ is three times differentiable, and  $(Sf)(x)<0$ whenever $f'(x)\neq 0$, where $Sf$ denotes the Schwarzian
derivative of $f$, defined by
  $$
(Sf)(x)=\frac{f'''(x)}{f'(x)}-\frac{3}{2}\left(\frac{f''(x)}{f'(x)}\right)^2.
$$
\end{enumerate}

The following proposition is a consequence of Singer's results \cite{si}.
\begin{proposition}
\label{psinger}
Assume that $g:[a,b]\to [a,b]$  satisfies (S) and $g'(x)<0$ for all $x\in[a,b]$. Let $K$ be the unique fixed point of
$g$ in $[a,b]$. Then,
\begin{itemize}
\item If $|g'(K)|\leq 1$, then $\lim_{n\to\infty}g^{n}(x)=K\, ,\;\forall\, x\in [a,b].$
\item If $|g'(K)|>1$, then there exists a globally attracting $2-$cycle $\{p,q\}$. More precisely, $g(p)=q$, $g(q)=p$,
$p\neq q$, and
$$\lim_{n\to\infty}g^{2n}(x)=\left\{
\begin{array}{l} p\, \text{ if }\;  x<K\\
\noalign{\medskip}   q\, \text{ if } \; x>K.\end{array}
\right.
$$
\end{itemize}
\end{proposition}

For a detailed proof of the first statement in a more general situation, see, e.g., \cite[Proposition 3.3]{ltt}. The
second statement can be easily proved using the same arguments.

\begin{remark}
\label{rr}
As noticed above, condition (S) holds for the Mackey-Glass and Nicholson's blowflies models, among others. See, for
example, \cite{gbt, gt, mpn2}. Thus, in these models, the interval $\bar J$ given in Lemma \ref{l1} reduces to $\{K\}$
if $|g'(K)|\leq 1$, whereas $\bar J=[p,q]$ if $|g'(K)|>1$.
\end{remark}

 The next result shows that  $\bar J$ is actually  the sharpest invariant and attracting interval containing the global
 attractor of (\ref{1}) for all $\tau$ if (L) and (S) hold. We emphasize that in the limit case $g^2(x_0)=x_0$ in (L),
 we have $\bar\alpha=\alpha=x_0$, $\bar\beta=\beta=g(x_0)$. Thus, as suggested in \cite{rw}, the intervals $\bar J$ and
  $J$ coincide in this special situation.

\begin{proposition}
\label{psharp}
Assume that (L) is fulfilled, and (S) holds in the interval $[\bar\alpha, \bar\beta]$.
If $|g'(K)|>1$ then, for any $\xi>0$, there exists a sufficiently large $\tau_{\xi}$ such that the interval
$[\bar\alpha+\xi, \bar\beta-\xi]$ is not attracting for (\ref{1}) if $\tau>\tau_{\xi}$.
\end{proposition}
\proof We make use of Theorems 2.1 and 2.2 in \cite{mpn2} (see also \cite{mpn}). Define  $\varepsilon=(\mu\tau)^{-1}$.
A direct application of \cite[Theorem 2.1]{mpn2}  provides an $\varepsilon_0>0$ such that for each
$\varepsilon< \varepsilon_0$, Eq. (\ref{1})
 possesses a slowly oscillating periodic solution $p(t)$ and there exist constants $\zeta,\omega>0$  satisfying $p(t)\in\bar J$,
 for all $t$, $p(0)=p(\zeta)=p(\omega)=K$, $p(t)>K$ in $(0,\zeta)$, $p(t)<K$ in $(\zeta,\omega)$, $p(t)=p(t+\omega)$ for
 all $t$.

 Next, \cite[Theorem 2.2]{mpn2} states that, given any $\delta>0$, there exist $\varepsilon_1>0$, $k>0$ such that
 $|p(t)-\bar\alpha|\leq\delta$ in $[\varepsilon k, \zeta-\varepsilon k]$ for all $\varepsilon< \varepsilon_1$. Choosing
 $\delta=\xi/2$, it is clear that the periodic solution $p(t)$ is not attracted by the interval
 $[\bar\alpha+\xi, \bar\beta-\xi]$ if $\tau>\tau_{\xi}=(\mu\varepsilon_1)^{-1}$.
\qed

\begin{remark}
What is important in the proof of Proposition \ref{psharp} is the fact that $g$ has a unique globally attracting
$2$-periodic solution. Condition (S) implies this fact, although it is not a necessary condition. For example, function
$$
g(x)=30 \frac{x + x^{5/2}}{2 + 35 x^3}
$$
satisfies (L) and (U), $|g'(K)|>1$, and it has a unique globally attracting $2$-cycle defined by $\bar\alpha=0.728449$,
$\bar\beta=2.2822$. The conclusion of Proposition \ref{psharp} holds although $(Sg)(x)>0$ for $x>2.02$. Notice that
$[\alpha,\beta]=[0.515162,3.62133]$.

On the other hand, if $g$ does not have a unique globally attracting $2$-periodic solution, under condition (L), interval
$\bar J$ is still the smallest globally attracting interval for the difference equation $x_{n+1}=g(x_n), \, n=0,1,\dots$
We conjecture that the conclusion of Proposition \ref{psharp} remains valid without assumption (S).
\end{remark}

As a consequence of Corollary \ref{c3}, Remark \ref{r1} and Propositions \ref{psinger} and \ref{psharp}, we get the
following dichotomy for Equation (\ref{1}) under conditions (L) and (S):
\begin{theorem}
\label{tdic}
Assume that (L) is fulfilled and (S) holds in the interval $[\alpha,\beta]$. Then exactly one of the following holds:
\begin{itemize}
\item[(1)] $|g'(K)|\leq 1$ and the global attractor of (\ref{1}) for all values of the delay $\tau$ is $\{K\}$.
\item[(2)] $|g'(K)|>1$ and the sharpest invariant and attracting interval containing the  global attractor of (\ref{1})
for all values of the delay $\tau$ is $[\bar{\alpha},\bar{\beta}],$ where $\{\bar{\alpha},\bar{\beta}\}$ is the unique
$2$-cycle of $g$ in $[\alpha,\beta]$.
\end{itemize}
\end{theorem}

According to Theorem \ref{trw}, when condition {\rm (L)}  does not
hold, it is still possible to find a delay-dependent condition {\rm
(L$_{\tau}$)} under which  every solution of (\ref{1}) enters the
domain where $f'$ is negative and the theory of monotone delayed
feedback can be applied to describe the asymptotic behaviour of our
equation with unimodal feedback \cite{rw}. Next we give a different
condition in the same direction. The proof is based on the following lemma proved in
\cite{gt} (see also \cite[Lemma 5.1]{ltt}).
\begin{lemma}
\label{lgt}
Assume that (U) holds, $g'(0)>1$, and $K>x_0$. Then,  for every nonnegative and nonzero solution
$x(t)$ of (\ref{1})  there exist
finite positive limits $$M = \limsup_{t \to \infty}x(t)\, ,\, m =
\liminf_{t \to \infty}x(t).$$
Moreover,  $[m,M] \subset g_1([m,M]),$ where 
\begin{equation}
\label{g1}
g_1(x):=\left(1-e^{-\mu\tau}\right) g(x)+e^{-\mu\tau}K.
\end{equation}
\end{lemma}

\begin{theorem}
\label{ddc}
Assume that the following
condition holds:
\begin{enumerate}
\item[{\rm (L'$_{\tau}$)}] $g_1^2(x_0)>x_0,\,$
where $g_1(x)=\left(1-e^{-\mu\tau}\right) g(x)+e^{-\mu\tau}K$.
\end{enumerate}

Then, every solution of (\ref{1}) enters the domain where $f'$ is negative.
\end{theorem}

\proof Let $x(t)$ be a nonnegative nonzero solution of (\ref{1}),  $m=\displaystyle\liminf_{t\to\infty}x(t)$ and  $M \displaystyle =\limsup_{t\to\infty}x(t)$. 
Notice that $g_1(x)$ lies between $g(x)$ and $K$. In particular, $g_1(x_0)>x_0$. On the other hand, since $f$ satisfies (U), it is clear that $g_1$ also meets the same condition, except that $g_1(0)>0$. Moreover, $g_1$
has the only fixed point $K$.  
Since, by  Lemma \ref{lgt}, $[m,M]\subset g_1([m,M])$, we obtain $m\leq K \leq M$, thus $m \leq g_1(m)$, $M \geq g_1(M)$
 and the following inequalities must hold:
$$
\min\{g_1(M),g_1^2(x_0)\}\leq m, \quad M\leq g_1(x_0).
$$
Since $g_1$ is decreasing in $[M,g_1(x_0)]$,  $M\leq g_1(x_0)$ implies  $g_1^2(x_0)\leq g_1(M)$ and then 
$$g_1^2(x_0)=\min\{g_1(M),g_1^2(x_0)\}\leq m.$$
 Finally,  {\rm (L'$_{\tau}$)} implies that $m\geq g_1^2(x_0)>x_0$, and therefore $x(t)>x_0$
for all sufficiently large $t$. The proof is complete.
\qed


As a byproduct of the proof of Theorem \ref{ddc}, we get the following result, which is of independent interest to obtain sharper bounds for the global attractor of Eq. (\ref{1}).
\begin{corollary}
\label{cornew}
Let $g_1$ be the function defined in (\ref{g1}). If (U) holds, $g'(0)>1$, and $K>x_0$, then the interval $[g_1^2(x_0),g_1(x_0)]$ contains the global attractor of (\ref{1}). 
\end{corollary}

\begin{figure}
\label{gprofile} \centering
\includegraphics[totalheight=2.5in]{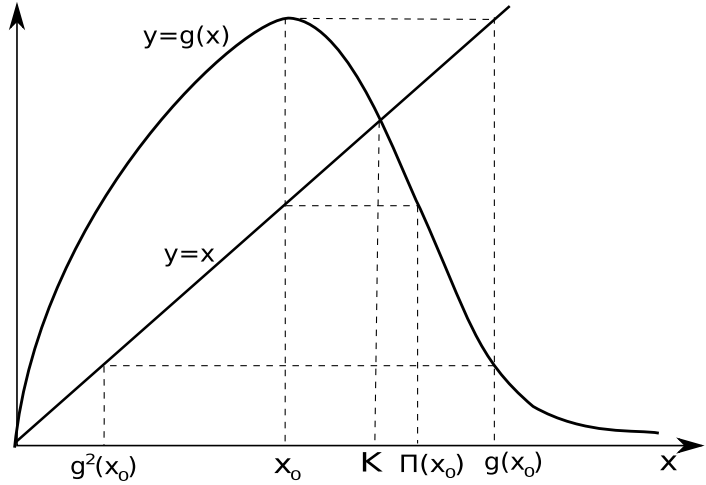}
\caption{The profile of the function $g$ in a situation where condition (L) does not hold.}
\end{figure}

Notice that {\rm (L'$_{\tau}$)} always holds if {\rm (L)} is satisfied, since $g_1^2(x_0)>\alpha=g^2(x_0)$. Next we show that when {\rm (L)} does not hold,  {\rm (L'$_{\tau}$)}  is sharper than  {\rm (L$_{\tau}$)}, and therefore Theorem \ref{ddc}  improves Theorem 3.8  of \cite{rw}.

\begin{proposition}
Assume that (L) does not hold. Then ($L_\tau$) implies ($L'_\tau$).
In other words, ($L'_\tau$) gives a better estimate than ($L_\tau$)
for the possible delays that still guarantee that every solution
enters the domain where $f'$ is negative.
\end{proposition}
\proof If (L) fails, then $g^2(x_0)<x_0$. Thus we have
$g^2(x_0)<x_0<K<\Pi(x_0)<g(x_0)$ (as depicted in Figure 1). Using
the notation $\theta=\mu\tau \geq 0$, ($L_\tau$) is equivalent with
$$x_0+\theta (g(x_0)-g^2(x_0)) < \Pi(x_0)$$
and ($L'_\tau$) can be written as
$$(1-e^{-\theta}) g[(1-e^{-\theta}) g(x_0)+e^{-\theta} K]+e^{-\theta} K > x_0.$$

First we suppose that $(1-e^{-\theta}) g(x_0)+e^{-\theta} K \leq
\Pi(x_0)$. Notice that $g(x_0)>(1-e^{-\theta}) g(x_0)+e^{-\theta} K
> K$, and $g$ is decreasing in in $[K,g(x_0)]$, so $$g[(1-e^{-\theta})
g(x_0)+e^{-\theta} K]\geq g(\Pi(x_0))=x_0,$$ and by
$$(1-e^{-\theta}) g[(1-e^{-\theta}) g(x_0)+e^{-\theta} K]+e^{-\theta} K \geq
(1-e^{-\theta}) x_0+e^{-\theta} K   > x_0,$$ ($L'_\tau$) holds. So
far we have not used $(L\tau)$, $(1-e^{-\theta}) g(x_0)+e^{-\theta}
K \leq \Pi(x_0)$ always implies $(L_\tau')$.

Now we consider the remaining case $(1-e^{-\theta})
g(x_0)+e^{-\theta} K > \Pi(x_0)$. Since ($L_\tau$) holds, we have
$(1-e^{-\theta}) g(x_0)+e^{-\theta} K
>x_0+\theta (g(x_0)-g^2(x_0))$. On the other hand,
$$g[(1-e^{-\theta}) g(x_0)+e^{-\theta} K] > g^2(x_0).$$ Therefore,
\begin{eqnarray*} (1-e^{-\theta})g[(1-e^{-\theta}) g(x_0)+e^{-\theta}
K]+e^{-\theta} K > (1-e^{-\theta}) g^2(x_0)+e^{-\theta} K  =  \\ =
(1-e^{-\theta}) g(x_0)+e^{-\theta} K+(1-e^{-\theta})(
g^2(x_0)-g(x_0)) > \\
>x_0+\theta (g(x_0)-g^2(x_0))+(1-e^{-\theta})( g^2(x_0)-g(x_0)) =\\
=x_0+ (g(x_0)-g^2(x_0))(\theta+e^{-\theta}-1)>x_0,\end{eqnarray*}
where in the last step we used that $g(x_0)-g^2(x_0)>0$ (this
follows because (L) does not hold), and that the function
$h(\theta)=\theta+e^{-\theta}-1$ is nonnegative for $\theta \geq 0$.
\qed

\begin{remark} Notice that {\rm (L$_{\tau}$)} always fails if
$\theta \geq 1$ (or equivalently $\tau \geq 1/\mu$), but we have not
used this fact in the proof.
\end{remark}

\section{Examples}
\setcounter{equation}{0}

In this section, we use the Nicholson's blowflies equation and the
Mackey-Glass equation with different parameters in order to
illustrate our results in Section 2.

After a change of variables, one can always write (\ref{nic}) in the form
\begin{equation}
\label{nic2}
x'(t)=- \mu x(t)+ x(t-\tau)e^{-x(t-\tau)}.
\end{equation}

It is well known that there is a unique equilibrium $x=0$ if $\mu\geq 1$, and it attracts all nonnegative solutions.
If $\mu<1$, there is a positive equilibrium $K=-\ln(\mu)$, and $x=0$ becomes unstable. Moreover, $K$ is globally
attracting for all values of the delay if $\mu\in(e^{-2},1)\approx (0.13533, 1)$. Next, using Theorem 9.3 in \cite{mpn2},
 it follows that condition (L) holds for $\mu\in(\nu, e^{-2})$, where $\nu\approx 0.10472$ is defined by the relation
  $\nu=e^{-1}v^{-1}$, $v$ being the unique solution greater than $2$ of equation $v^2=e^{v-1}.$

For $\mu=0.13$, condition (L) holds, and the  invariant and attracting interval for (\ref{nic2}) given by Theorem 3.5
in \cite{rw} is
$$
[\alpha,\beta]=[g^2(x_0),g(x_0)]=[1.2848,2.8298].
$$
This interval may be strengthened, since the smallest invariant and attracting interval for (\ref{nic2}) independent of
the delay $\tau$ is given by the unique $2$-cycle of $g$, that is,
$$
[\bar\alpha, \bar\beta] \approx [1.54796, 2.53248].
$$
According to Proposition \ref{psharp}, there exists a slowly
oscillating periodic solution of (\ref{nic2}) whose minimum and
maximum values get closer and closer to $\bar\alpha$ and $\bar\beta$
as $\tau$ tends to infinity. See Figure 2.A, where two distinct
solutions are presented and the horizontal lines indicate
$\bar\alpha$ and $\bar\beta$.

\begin{figure}
[hb] 
\label{fignic1} \centering
\includegraphics[totalheight=2.1in]{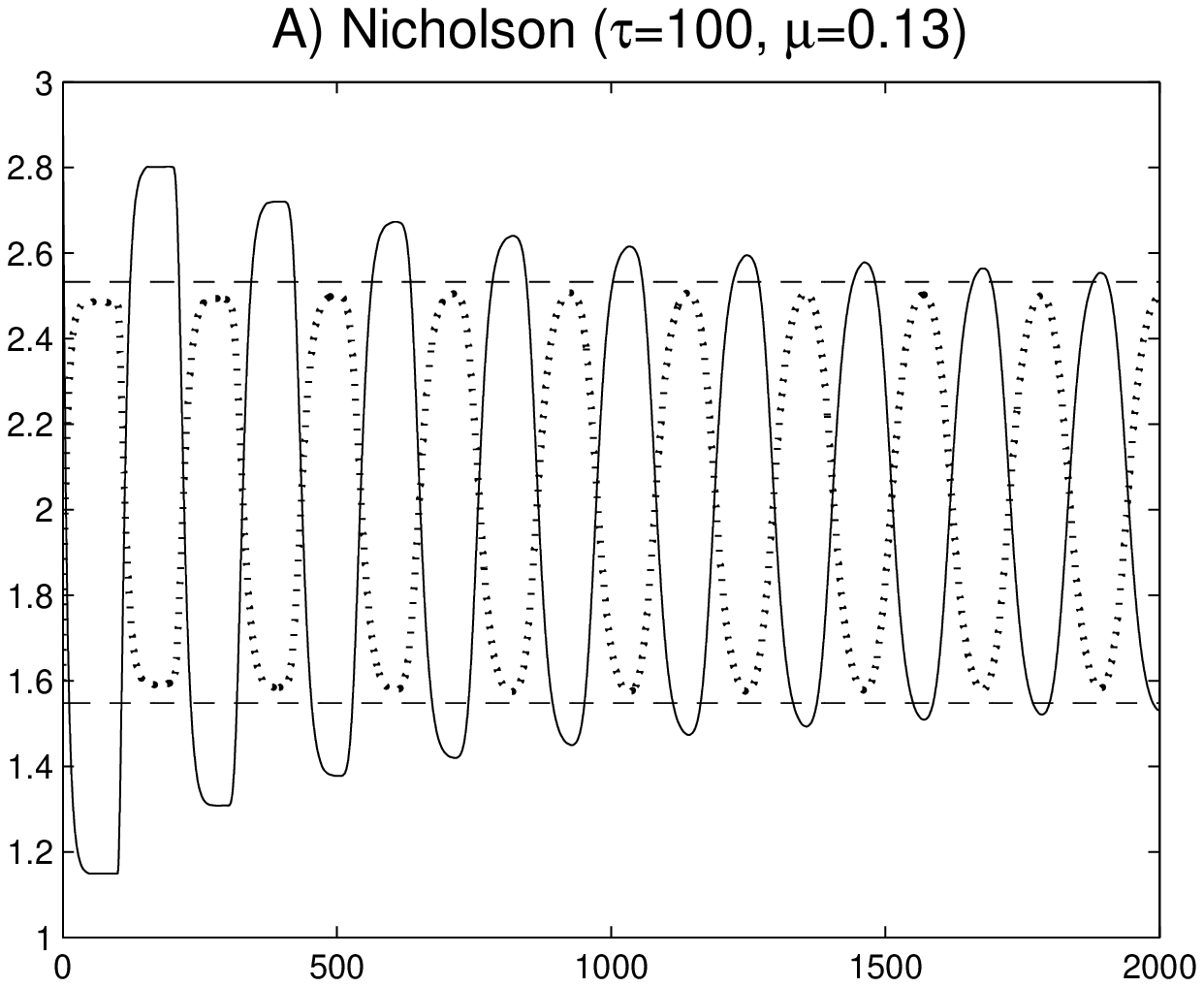}
\includegraphics[totalheight=2.1in]{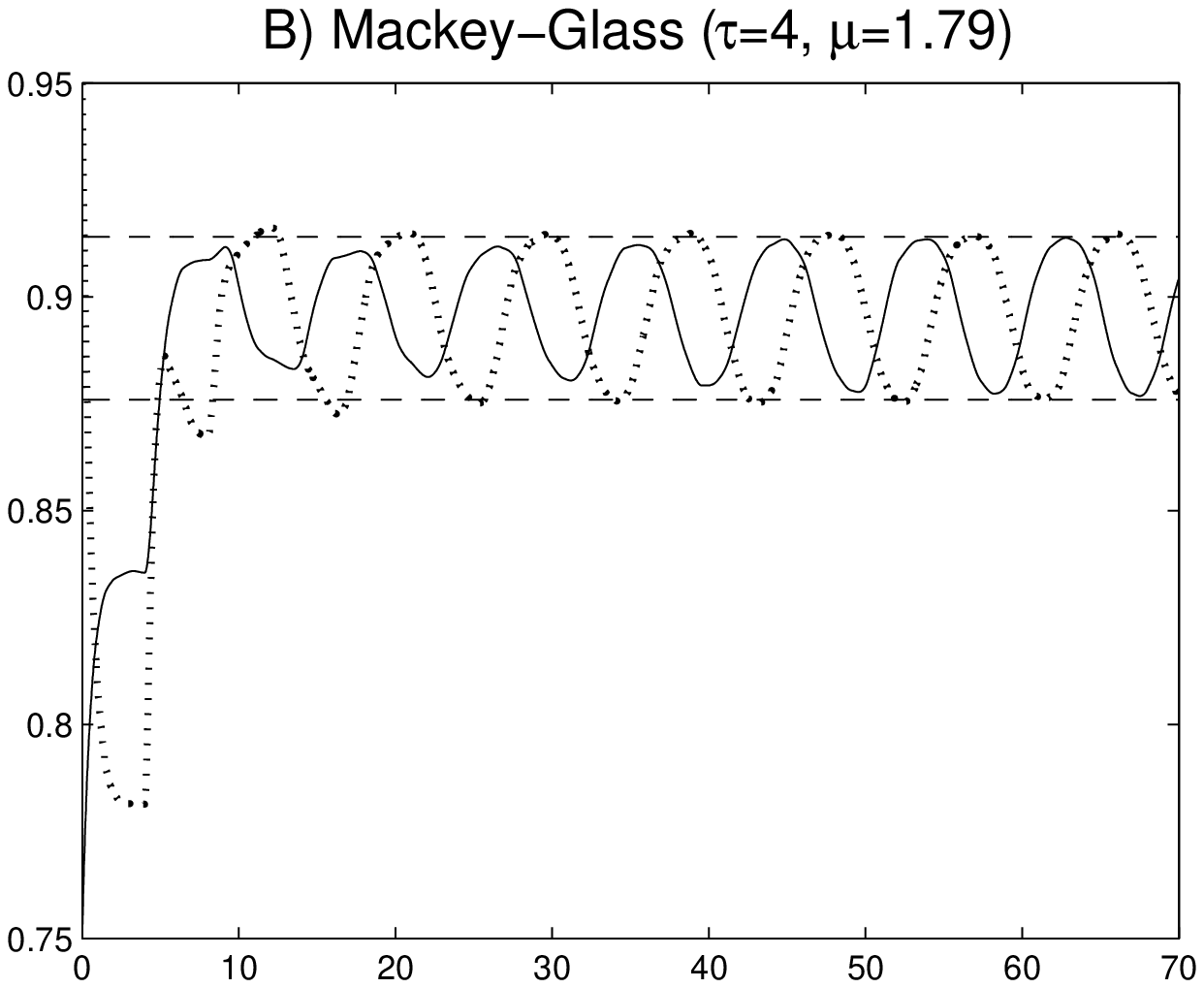}
\caption{Nicholson's blowflies equation and Mackey-Glass equation
when the condition (L) holds}
\end{figure}

For $\mu=1/16=0.0625$, condition (L) does not hold (indeed,   $g^2(x_0)=0.2616<x_0=1).$ Since $\Pi(x_0)=4.21007$, we get
$$
\tau^*=\frac{\Pi(x_0)-x_0}{\mu(g(x_0)-g^2(x_0))}=0.570734.
$$

On the other hand, one can check that condition {\rm (L'$_{\tau}$)}
in Theorem \ref{ddc} holds for $\tau<\tau_*=1.46534$, showing how {\rm
(L'$_{\tau}$)} gives an estimate significantly sharper than   {\rm (L$_{\tau}$)}.

It is interesting to notice that from Theorem 2.1 in \cite{ltt} it
follows  that the positive equilibrium $K$ is globally attracting
for (\ref{nic2}) if $\tau<1.1935$. Hence, the information provided
by Theorem \ref{trw} is not very useful here, whereas Theorem
\ref{ddc} can be applied for $\tau$ between $1.1935$ and $1.46534$.

\begin{figure}
\label{figmg1} \centering
\includegraphics[totalheight=2.1in]{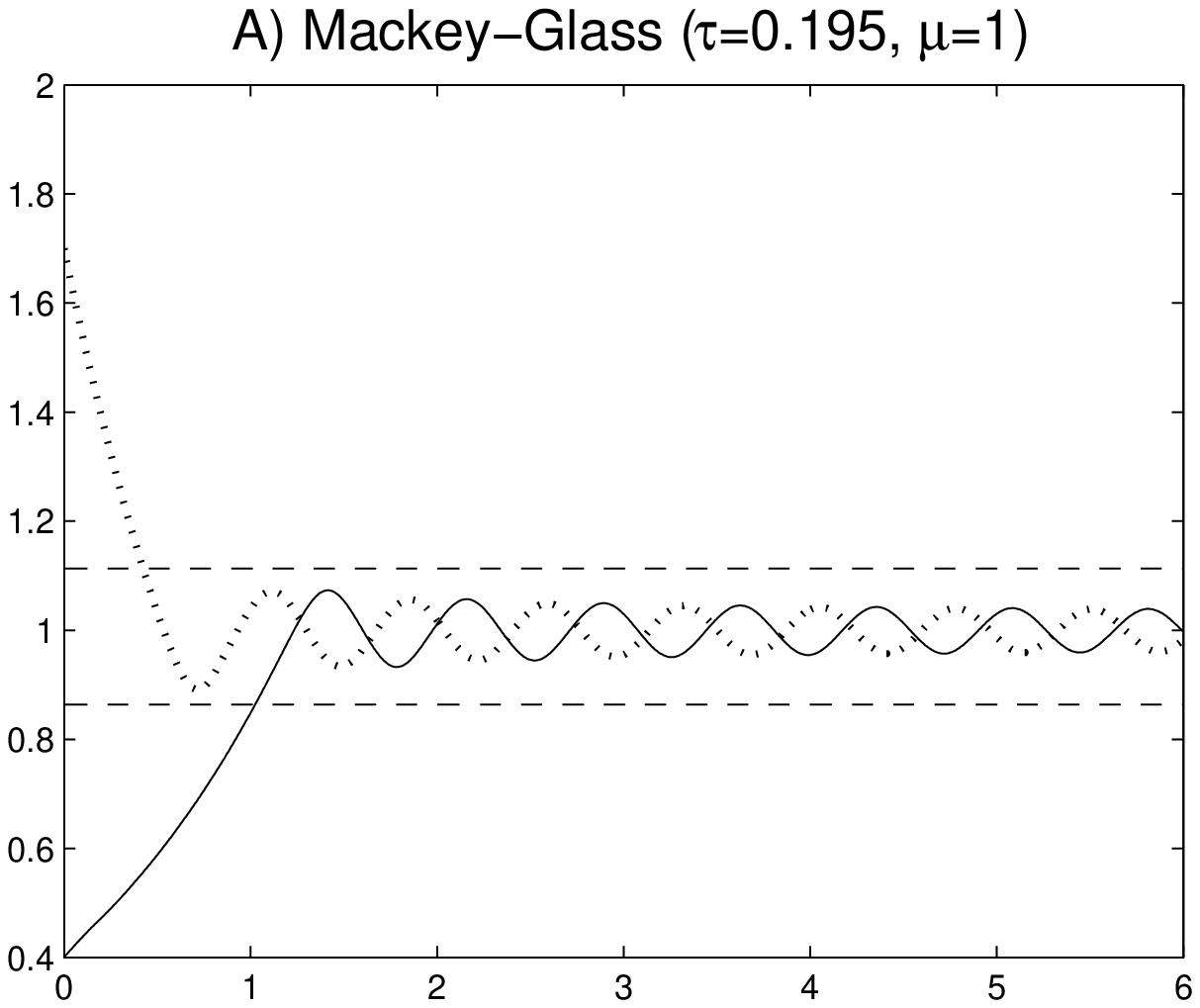}
\includegraphics[totalheight=2.1in]{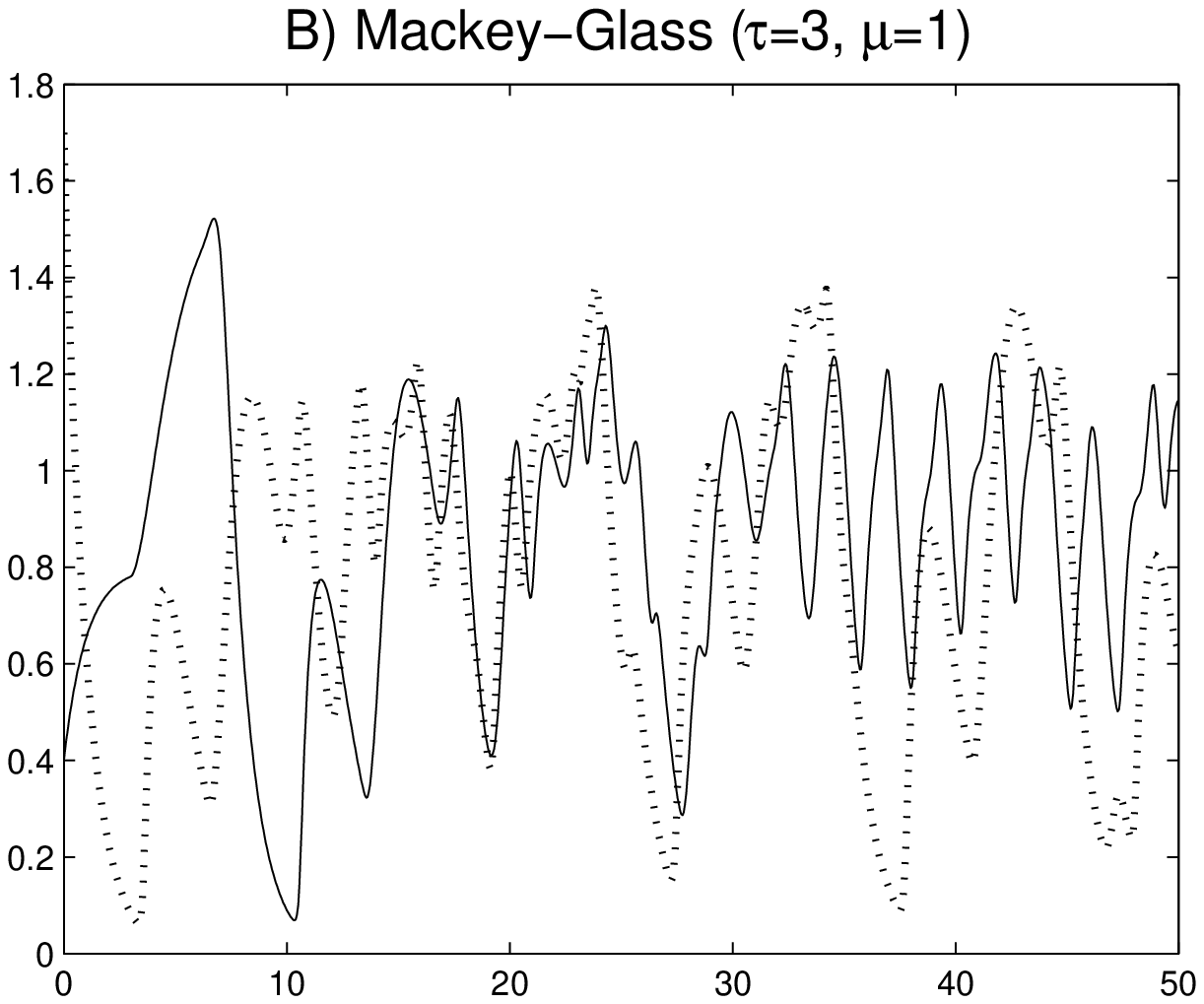}
\caption{Illustration of Theorem \ref{ddc}: Mackey-Glass equation with
different delays}
\end{figure}

Our next example is the Mackey-Glass equation

\begin{equation}
\label{mg1} x'(t)=-\mu x(t)+\frac{2x(t-\tau)}{1+x(t-\tau)^{20}}.
\end{equation}

The function $f(x)=2x(1+x^{20})^{-1}$ satisfies the unimodal
condition (U) with $x_0=0.863$. For $\mu\geq 2$,  there is a unique
equilibrium $x=0$ which attracts all nonnegative solutions. If
$\mu<2$, there is also a positive equilibrium
$K=(-1+2\mu^{-1})^{1/20}$. Next, $K\leq x_0$ for $\mu\in [1.9,2)$
(and hence $K$ is globally attracting), while $K>x_0$ for $\mu<1.9$.

Denote, as usual, $g=\mu^{-1}\, f$. One can check that $(L)$ holds if  $\mu> 1.774$ so in the interval
 $(1.774, 1.9)$ the dichotomy stated in Theorem \ref{tdic} applies. Since $|g'(K)|\leq 1$ if and only
 if $\mu\geq 1.8$, the equilibrium $K$ attracts all positive solutions for $\mu\in [1.8,1.9)$, whereas
  for each $\mu\in (1.774,1.8)$ we can determine the sharpest delay-independent interval
  $[\bar\alpha,\bar\beta]$ containing the global attractor of (\ref{mg1}) by finding the unique $2$-cycle
  of $g$ in the interval $[\alpha,\beta]=[g^2(x_0), g(x_0)]$. For example, setting $\mu=1.79$, we have
  $K=0.898$, $[\alpha,\beta]=[0.872,0.916]$ and the interval
$[\bar\alpha,\bar\beta]=[0.876,0.914]$ contains the global attractor
of (\ref{mg1}) for all values of the delay. See Figure 2 B, where
the horizontal lines represent $\bar \alpha$ and $\bar \beta$.

When (L) does not hold, we still can use Theorem  \ref{ddc}. For example, for $\mu=1$ the positive equilibrium
  $K=1$ loses its asymptotic stability for $\tau>0.188$. One can check that
condition {\rm (L'$_{\tau}$)} in Theorem~\ref{ddc} holds for
$\tau<\tau_*=0.195$, so for $\tau\in (0, \tau_*)$ every solution of
(\ref{mg1}) enters the domain where $f'$ is negative, while {\rm
(L$_{\tau}$)} is satisfied only for $\tau<\tau^*=0.092$. In this
situation apparently chaotic behavior can be observed for large
enough delays (see Figure 3 B with $\tau=3$), but Theorem~\ref{ddc}
and the results of \cite{rw} guarantee that complicated behavior is
not possible if $\tau<\tau_*=0.195$, as illustrated in Figure 3 A.
Moreover, we get a good bound for the global attractor of  (\ref{mg1}) from Corollary \ref{cornew}. Indeed, in this case the interval $[\alpha,\beta]=[0.00016,1.639]$ is improved up to $[g_1^2(x_0), g_1(x_0) ]=[0.864, 1.113]$. See Figure 3 A, where the horizontal lines indicate $g_1^2(x_0)$ and $g_1(x_0)$.

\section*{Acknowledgements}
E. Liz was partially supported by MEC (Spain)  and FEDER,  grant
MTM2007-60679. G. R\"{o}st was partially supported by the
Hungarian Foundation for Scientific Research, grant T 049516, NSERC
Canada and MITACS.

 \end{document}